\documentclass[11pt]{amsart}
\usepackage[foot]{amsaddr}
\usepackage{geometry}
\usepackage[T1]{fontenc}
\usepackage[english]{babel}
\addto\extrasenglish{
  
}
\usepackage{amssymb}
\usepackage{amsmath}
\usepackage{amsthm}
\usepackage{amsfonts}
\usepackage[shortlabels]{enumitem}
\usepackage{graphicx}
\usepackage{mathabx}
\usepackage{eurosym}
\usepackage{dsfont}
\usepackage{csquotes}
\usepackage{xparse}
\usepackage{accents}
\usepackage{color}
\usepackage{xcolor}
\usepackage{blindtext}
\usepackage{thmtools}
\usepackage{mathrsfs}
\usepackage{stmaryrd}
\usepackage{thm-restate}
\usepackage[
    backend=biber,
    style=alphabetic,
  ]{biblatex}

\usepackage{hyperref}
\usepackage{xurl}
\hypersetup{breaklinks=true}
\hypersetup{
	colorlinks = true,
	linkcolor = {purple},
	urlcolor = {blue},
	citecolor = {brown}
}
\hypersetup{
colorlinks = true,
linkcolor = {purple},
urlcolor = {blue},
citecolor = {brown}
}

\newtheorem{lem}{Lemma}

\newtheorem{cor}[lem]{Corollary}
\newtheorem{prop}[lem]{Proposition}

\newtheorem{remark}{Remark}

\newcommand{\RR}{{\mathbb R}}

\newcommand{\Pf}{\mathfrak{P}}
\newcommand{\XX}{\mathbb{X}}

\newcommand{\R}{\mathcal{R}}

\def\beq{\begin{equation}}
\def\eeq{\end{equation}}
\let\origautoref\autoref
\def\autoref#1{\textbf{\origautoref{#1}}}

\title{On the contraction properties of a pseudo-Hilbert projective metric}
\author{Maxime Ligonnière}\address{Institut Denis Poisson UMR 7013, Université de Tours, Université d’Orléans, CNRS France}
\email{maxime.ligonniere@univ-tours.fr}
\address{CMAP, CNRS, INRIA, École polytechnique, Institut Polytechnique de Paris, 91120 Palaiseau, France}

\begin{document}
\begin{abstract}
	In this note, we define a bounded variant of the Hilbert projective metric on an infinite dimensional space $E$ and study the contraction properties of the projective maps associated with positive linear operators on $E$. More precisely, we prove that any positive linear operator acts projectively as a $1$-Lipschitz map relatively to this metric. We also show that for a positive linear operator, strict projective contraction is equivalent to a property called uniform positivity. 
\end{abstract}
\maketitle
\medskip 
Let $\XX$ be a set of arbitrary cardinality and consider a vector space $E\subset\RR^\XX$. The projective space $\Pf$ associated with $E$ is defined as the set of equivalence classes $\Pf=E/\R$, where $\R$ is the equivalence relation such that for any $f,g\in E$
\[ f\R g \Leftrightarrow \exists b\in \RR_+^*, f=b g\]
Let $\Pi:E\longrightarrow \Pf$ be the canonical projection. 
The Hilbert metric $d_H$, defined for example in \cite{busemann_projective_1953}, is a distance on the projective image $\Pi(C)$ of the positive cone $C=E\cap \RR_+^\XX-\{0\}$. Any linear map $M$ on $E$ which is positive, in the sense that $M(C)\subset C$, generates a projective action on $\Pi(C)$. We say that the projective action $M:\Pi(C)\longrightarrow \Pi(C)$ associated with a linear map operator $M$ is $k$-contracting, with $k<1,$ when it is $k$-Lipschitz with respect to $d_H$. When $M$ is $k$-contracting for some $k<1$, we say that it is contracting, or strictly contracting. As proved in \cite{birkhoff_extensions_1957}, the projective action on $\Pi(C)$ of any positive bounded linear operator is $1$-Lipschitz with respect to $d_H$.
We say that a positive, linear, bounded operator $M$ on $E$ is $A$-uniformly positive for some $A>1$ when there exists $h\in E \cap (\RR_+)^\XX$ such that, for any $f\in E$, there exists $ b(f)\geq 0$ satisfying
\begin{equation}
A^{-1} b(f) h\leq Mf \leq A b(f).
\end{equation}
Birkhoff \cite{birkhoff_extensions_1957} shows that uniformly positive operators are contracting with respect to the Hilbert metric. This is useful to prove the existence a fixed point for the projective action of such an operator, that is an eigenvector for its linear action. It also allows to study the ergodicity properties of semi-groups of uniformly positive operators.

In this note, we focus on a bounded variant $d$ of this Hilbert distance, introduced in a finite dimensional setup in \cite{hennion_limit_1997}. As shown in \cite{hennion_limit_1997}, any matrix with non negative coefficients is $1$-Lipschitz with respect to $d$, and an explicit formula allows to compute the contraction rate of a matrix in terms of its coefficients is provided. This yields a sufficient condition for a matrix to be strictly contracting with respect to 
$d$. In this note, we provide an elementary construction of this pseudo Hilbert distance $d$ in any dimension, as well as a study of some of its properties. In particular, we prove in Proposition \ref{prop:contraction} that all linear, positive and bounded operators are $1$-Lipschitz and in Proposition \ref{prop:strict contraction} that the contracting operators with respect to $d$ are exactly the uniformly positive ones. To the best of our knowledge, it was not yet proven in the literature that uniform positivity is not only sufficient but also necessary for an operator to be strictly contracting with respect either to the pseudo-Hilbert or the Hilbert metric.

\section{Defining a pseudo-Hilbert distance}
Consider the partial order $\leq$ defined by $f\leq g \Leftrightarrow \forall x\in \XX, f(x)\leq g(x)$.

For any $f,g\in C$, we define 
\[ \aleph(f,g)=\sup\{ b \geq 0 | b f \leq g \} \in [0,\infty]=\inf\left\{ \frac{g(x)}{f(x)}, x\in \XX, f(x)\neq 0\right\}, \]
\[m(f,g)=\aleph(f,g)\aleph(g,f).\]
Notice that the set $\{ b \geq 0 | b f \leq g \}$ clearly is a sub-interval of $\RR_+$, which contains $0$.
Moreover, it holds 
\begin{lem} \label{lem:ratios}
For any $f,g,h\in C$, any $\alpha,\beta>0$
\begin{enumerate}
    \item[i)] $\aleph(f,g)<\infty$ and $\{ b \geq 0 | b f \leq g \}=[0,\aleph(f,g)].$ 
    \item[ii)] $m(f,g)=m(g,f)$
    \item[iii)] $m(\alpha f, \beta g)= m(f,g)$,
    \item[iv)] $m(f,g)m(g,h)\leq m(f,h).$
    \item[v)] $0\leq m(f,g)\leq 1.$
    \item[vi)] $m(f,g)=1\Leftrightarrow f \R g$. If $m(f,g)=1$, then $f=\aleph(g,f)g$.
    
\end{enumerate}
\end{lem}
\begin{proof}
\begin{enumerate}
\item[\textit{i)}]  Since $0\notin C$, $f\neq 0$, thus, there exists $x\in \XX$ such that $f(x)>0$. For $b$ large enough, $b f(x)>g(x)$ which prevents $b f \leq g$. Since $\{ b \geq 0 | b f \leq g \}$ is an interval, this implies that it is bounded, i.e. $\aleph(f,g)<\infty.$ Thus $b f \leq g$, for any $b\in \left[0,\aleph(f,g)\right)$. This yields $\aleph(f,g) f \leq g$, thus \[\aleph(f,g)\in\{ b \geq 0 | b f \leq g \}=[0,\aleph(f,g)].\]
\item[\textit{ii)}] This symmetry property is straightforward from the definition of $m(f,g).$
\item[\textit{iii)}] This derives directly from \[\aleph(\alpha f,\beta g)= \frac{\beta}{\alpha} \aleph(f,g).\]
\item[\textit{iv)}] Combining $\aleph(f,g) f \leq g$, and $\aleph(g,h) g \leq h $, we get 
\[ \aleph(g,h)\aleph(f,g) f \leq \aleph(g,h) g \leq h,\]
thus $\aleph(g,h)\aleph(f,g)\in \{b \geq 0 | b f\leq h\}=[0,\aleph(f,h)]$. This yields $\aleph(g,h)\aleph(f,g)\leq \aleph(f,h).$
For a similar reason, it holds $\aleph(g,f)\aleph(h,g)\leq \aleph(h,f).$
Multiplying these two inequalities yields $m(f,g)m(g,h)\leq m(f,h)$.
\item[\textit{v)}] Notice that taking $f=h$ in \textit{iv)} yields 
\[m(f,g)m(g,f)=m(f,g)^2\leq m(f,f).\]
Since $m(f,f)=1$, then $m(f,g)^2\leq 1$, which implies that $m(f,g)\leq 1$.
\item[\textit{vi)}] From \textit{iii)}, it is clear that if $f\R g$ then $m(f,g)=m(f,f)=1$. 
Suppose now that $m(f,g)=1$. This implies that $\aleph(f,g)\neq 0, \aleph(g,f)\neq 0$ and $\aleph(f,g)^{-1}=\aleph(g,f)$. Moreover, it holds
\[ \aleph(f,g) f \leq g\]
and \[ \aleph(g,f)g\leq f.\]
Thus 
\[ \aleph(f,g)\aleph(g,f) f = f \leq \aleph(g,f)g \leq f.\]
Therefore, $f=\aleph(g,f) g$.
\end{enumerate}
\end{proof}

Let $\bar{f},\bar{g}\in \Pi(C)$ be two halflines and $f\in \bar{f}$, $g\in\bar{g}$ two points of those halflines, we set 
\[d(\bar{f},\bar{g})=\frac{1-m(f,g)}{1+m(f,g)}.\]

\begin{prop}
$d$ is a well defined distance on $\Pi(C)$. It is bounded by $1$.
\end{prop}
\begin{proof}
Assertion \textit{iii)} of Lemma \ref{lem:ratios} implies that $m(f,g)$ only depends on the equivalence classes $\Pi(f), \Pi(g)$ and not on the choice of $f$ and 
$g$ inside those classes. Moreover since $m(f,g)\in [0,1]$, and $\phi : [0,1]\longrightarrow [0,1], s \mapsto \frac{1-s}{1+s}$ is well defined, then $d$ is well defined. The map $d$ is clearly symmetric, nonnegative and bounded by $1$.
Since the map $\phi$ is strictly decreasing, continuous, $\phi(0)=1$ and $\phi(1)=0$, we obtain 
\[d(\bar{f},\bar g )=0\Leftrightarrow m(f,g)=1 \Leftrightarrow f\R g \Leftrightarrow \bar{f}=\bar{g}.\]
Moreover, it can be checked that $\phi(st)\leq \phi(s)+\phi(t)$, for any $s,t\in [0,1]$. Combining this with the fact that $\phi$ is decreasing and point \textit{iv)} of Lemma \ref{lem:ratios} yields the triangular inequality. 
\end{proof}

\begin{remark}
In many references such as \cite{birkhoff_extensions_1957}, one rather considers the Hilbert metric $d_H(\bar{f},\bar g )=\vert \log m(f,g) \vert$. We prefer using the pseudo Hilbert metric because its boundedness makes it more convenient.
\end{remark}
\begin{prop} \label{prop:segments}
For any $\bar{f}\neq\bar{g}\in \Pi(C)$, the following claims hold
\begin{itemize}
    \item there exists  $f\in \bar{f},g\in \bar{g}$ such that the intersection of the line $(f,g)$ with the cone $C$ is reduced to a line segment $[u,v]$, with $u\neq v$.
Fix now such points $f,g,u,v$.
    \item For any $h\in C$, if $h$ is coplanar with $u$ and $v$ (or equivalently with $f$ and $g$), then the line segment $[u,v]$ intersects the vector line $(0,h)$ in a single point of the cone.
    \item For any point $h\neq 0$ coplanar with $u$ and $v$, with coordinates $h=h_1u+h_2v$ then $h\in C$ iff $h_1\geq 0$ and $h_2\geq 0$.
\end{itemize} 
\end{prop}
\begin{proof}
Let $f\in\bar{f}$ and $g\in\bar{g}$. For any $t\in\RR$, we note \[h_t=tf+(1-t)g\in E.\]
Consider $S=\{t\in \RR | h_t\in C\}$. $S$ is clearly an non empty interval in $\RR$ since $0,1\in S$ and $C$ is convex. In order to prove that $S$ is a segment of line, it is enough to show that there exists $x,x'\in\XX$ such that $f(x)<g(x)$ and $g(x')<f(x')$. Indeed, if this is the case, since $h_t=g+t(f-g)$, then \[h_t(x)\underset{t\rightarrow +\infty}{\longrightarrow -\infty} \text{ and } h_t(x')\underset{t\rightarrow -\infty}{\longrightarrow +\infty}.\]
Let us now prove that we can choose $f\in \bar{f}$ and $g\in \bar{g}$ such that there exists $x,x'\in\XX$ satisfying $f(x)<g(x)$ and $f(x')>g(x')$. Consider two points $\tilde{f}\in \bar{f}, \tilde{g}\in \bar{g}$. Since $\bar{f}\neq \bar{g}$, in particular $\tilde{f}\neq \tilde{g}$, thus there exists $x\in\XX$ such that $\tilde{f}(x)\neq \tilde{g}(x)$. Then in particular, $\tilde{f}(x)$ or $\tilde{g}(x)$ is nonzero. Without loss of generality, suppose $\tilde{g}(x)\neq 0$. Then, since $\bar{f}\neq \bar{g}$, there exists $x'\in\XX$ such that \[\tilde{f}(x')\neq \frac{\tilde{f}(x)}{\tilde{g}(x)} \tilde{g}(x').\]
If $\tilde{g}(x')=0$, then for any $b>0$, $b\tilde{f}(x')>0=\tilde{g}(x')$, and for $b >0$ small enough, $b \tilde{f} (x)< \tilde{g}(x)$, thus setting $f=b \tilde{f}$ and $g= \tilde{g}$ suffices. 
In the contrary, if $\tilde{g}(x')\neq 0$, then,
\[ \frac{\tilde{f}(x')}{\tilde{g}(x')}\neq \frac{\tilde{f}(x)}{\tilde{g}(x)}.\]
Once again, without loss of generality, suppose that
\[\frac{\tilde{f}(x')}{\tilde{g}(x')}>\frac{\tilde{f}(x)}{\tilde{g}(x)}. \]
Then there exists $b>0$ such that 
\[\frac{b \tilde{f}(x')}{\tilde{g}(x')}>1> \frac{b \tilde{f}(x)}{\tilde{g}(x)}.\]
Then, setting $f=b \tilde{f}$, $g=\tilde{g}$, we get ${f}(x)<{g}(x)$ and ${f}(x')>{g}(x')$, thus $S$ is bounded and $(f,g)\cap C$ is a nonempty line segment $[u,v]$. Moreover, $f\neq g$ and $f,g\in[u,v]$, thus $u\neq v$.
Finally, let $h\in C$, suppose that $h$ is coplanar with $f,g$. Note that all points of the vector line $(0,h)$ have either only nonnegative coordinates or only nonpositive coordinates. Thus $f-g$ does not belong to this line. This proves that the lines $(f,g)$ and $(0,h)$ are not parallel. Since they are coplanar, they are secant. Moreover, their intersection point $q=f+t(g-f)$ has either all nonnegative (in which case this point is in $[u,v]\subset C$), or all nonpositive coordinates. If $t=0$ then $q\in C$, thus in $q\in[u,v]$. If $t>0$, since $f(x)<g(x)$, then $q$ as a positive coordinate $q(x)$. Similarly, if $t<0$, the coordinate $q(x')$ is positive. Thus if $t\neq 0$, then the coordinates of $q$ are not all nonpositive, consequently they are all nonnegative and $q\in C\cap(f,g)=[u,v].$
Let now $h\neq 0$ be a point coplanar with $u$ and $v$. Since $u,v\in C$, if $h_1,h_2\geq 0$, then $h(x)=h_1u(x)+h_2v(x)\geq 0$ for any $x\in \XX$, thus $h\in C$. Conversely, if $h\in C$ the semi-line $\bar{h}$ intersects $[u,v]$. Thus there exists a point $\tilde{h}\in \bar{h}\cap [u,v]$. This point $\tilde{h}$ has non negative coordinates in the basis $(u,v)$ of the plan. Since $h$ is positively colinear with $\tilde{h}$, then $h_1,h_2\geq 0$.
\end{proof}
\begin{prop}\label{prop:expression_distance}
Let $\bar{f},\bar{g}\in \Pi(C)$. Consider a line $\Delta\in E$ that intersects both $\bar{f}$ and $\bar{g}$ and suppose that $\Delta\cap C$ is a line segment $[u,v]$, with $u\neq v$. Then any $f\in \bar{f}$ and $g\in \bar{g}$ are coplanar with $u$ and $v$. Moreover, noting $f=f_1 u + f_2 v$, $g=g_1 u + g_2 v$, it holds :
\begin{equation} \label{eq: expression distance} d(\bar{f},\bar{g})=\left|\frac{f_1g_2-f_2g_1}{f_1g_2+f_2g_1}\right|,\end{equation}
with the convention $0/0=0$.
\end{prop}
\begin{remark}
Applying the triangular inequality in \eqref{eq: expression distance} yields $d(\bar{f},\bar{g})\leq 1$ for any $\bar{f},\bar{g}$. Moreover, $d(\bar{f},\bar{g})=1$ iff exactly one of the quantities $g_1f_2$, $f_1g_2$ equals $0$. This implies in particular that $\bar{f}$ or $\bar{g}$ is contained in the edge of $C$.
\end{remark}
\begin{proof}
Since $\Delta=(u,v)$, then $\Delta$ is included in the vector plan spanned by $u,v$. The line $\Delta$ intersects the vector semi-lines $\bar{f},\bar{g}$, thus they are included in the plan spanned by $u,v$.
Note that the quantity 
\[\left|\frac{f_1g_2-f_2g_1}{f_1g_2+f_2g_1}\right|\]
is invariant by $(f,g)\mapsto (\alpha f,\beta g)$, for any positive real numbers $\alpha, \beta$, and by exchanging the coordinates $f_1$ with $f_2$ and $g_1$ with $g_2$. Therefore, it suffices to show that this equality is satisfied for a well chosen pair $(f,g)$ of elements of $\bar{f},\bar{g}$. Thus, we suppose in the sequel that $f,g$ are the intersection points of $ \Delta$ with $\bar{f},\bar{g}$. This implies that $f,g\in [u,v]$. Without loss of generality, we moreover assume that $u,f,g,v$ are in this order, in the sense that $f_1\geq g_1$, thus $g_2=1-g_1\geq 1-f_1=f_2$. It remains to prove that 
\[d(\bar{f},\bar{g})=\left|\frac{f_1g_2-f_2g_1}{f_1g_2+f_2g_1}\right|=\frac{f_1g_2-f_2g_1}{f_1g_2+f_2g_1}\geq 0.\]
This holds in the case where $f_1=g_1=0$ or in the case where $f_2=g_2=0$, with the convention $0/0=0$. Thus, since $f_1\geq g_1$ and $f_2\leq g_2$, it remains to study the case $f_1,g_2>0$. For any $b>0$, the last point of Proposition \ref{prop:segments} yields \[g\leq b f \Leftrightarrow g-b f =(g_1-b f_1) u +(g_2-b f_2)v\in C\Leftrightarrow g_1\geq b f_1\text{ and }g_2\geq b f_2.\]
Since $f_1\geq g_1$ and $f_2\leq g_2,$ it holds $f_1g_2\geq f_2g_1$ and thus 
\[ \frac{f_1}{g_1} \geq \frac{f_2}{g_2},\]
with the convention $f_1/g_1=\infty$  if $g_1=0$.
Thus, it holds \[ g\leq b f \Leftrightarrow b \leq \frac{f_2}{g_2},\]
and 
\[\aleph(f,g)=\frac{f_2}{g_2}\in(0,+\infty).\]
Similarly, \[\aleph(g,f)=\frac{g_1}{f_1}\in(0,+\infty).\]
Finally, it holds \[d(\bar{f},\bar{g})=\frac{1-\frac{g_1f_2}{f_1f_2}}{1+\frac{g_1f_2}{f_1f_2}}=\frac{{f_1f_2}-{g_1f_2}}{{f_1f_2}+{g_1f_2}}.\]
\end{proof}
\section{Contraction properties of positive linear mappings}
Let $M$ be a linear mapping on $E$, we assume that $M$ is positive, in the sense that $M(C)\subset C$. Then $M$ induces a projective action on $\Pi(C)$, defined by $M\cdot \Pi(f)= \Pi(Mf),$ for any  $f\in C$.
\begin{prop}\label{prop:contraction}
Let \[ c(M)=\sup \left\{d(M\cdot \bar{f}, M\cdot \bar{g} ),\bar{f},\bar{g}\in \Pi(C)\right\}\leq 1.\] Then for all $\bar{f},\bar{g}\in \Pi(C),$ it holds \begin{equation}\label{eq:lipschitz} d(M\cdot \bar{f}, M \cdot \bar{g})\leq c(M)d(\bar{f},\bar{g}).\end{equation}
\end{prop}
\begin{proof}
Let $\bar{f},\bar{g}\in C$. Suppose $M\cdot \bar{f} \neq M\cdot \bar{g}$, otherwise, inequality \ref{eq:lipschitz} is obvious. Let $f\in \bar{f},g \in \bar{g}$, $x\in M\cdot \bar{f}$ and $y\in M\cdot \bar{g}$ such that $(f,g)\cap C$ and $(x,y)\cap C$ are line segments. Let $a,b$ and $a_1,b_1$ be the respective extreme points of $(f,g)\cap C$, $(x,y)\cap C$. Notice that $(a,b)$ and $(a_1,b_1)$ are respectively the bases of two vector planes between which $M$ acts as an linear isomorphism. Let us decompose, $f$, $g$, $Ma$ and $Mb$ as follows :
$$f=f_1a+f_2b$$
$$g=g_1 a + g_2 b $$
$$Ma = \alpha a_1 + \gamma b_1$$
$$ Mb = \beta a_1 + \delta b_1.$$
Then
$$Mf=(f_1\alpha + f_2 \beta) a_1+(f_1 \gamma + f_2 \delta )b_1$$
$$Mg=(g_1\alpha + g_2 \beta)a_1+(g_1 \gamma + g_2 \delta) b_1.$$
The line $(a_1,b_1)$ intersects both $M\cdot \bar{f}$ and $M \cdot \bar{g}$, and $(a_1,b_1)\cap C= [a_1,b_1]$. Moreover, $Mf\in M\cdot \bar{f}$ and $Mg\in M\cdot \bar{g}$. Thus, by Proposition \ref{prop:expression_distance}, it holds 
\[d(M\cdot \bar{f}, M \cdot \bar{g})=\left|\frac{(f_1\alpha + f_2 \beta)(g_1 \gamma + g_2 \delta)-(f_1 \gamma + f_2 \delta )(g_1\alpha + g_2 \beta)}{(f_1\alpha + f_2 \beta)(g_1 \gamma + g_2 \delta)+(f_1 \gamma + f_2 \delta )(g_1\alpha + g_2 \beta)}\right|.\]
Simple calculations yields
\[d(M\cdot \bar{f}, M \cdot \bar{g})=\left|\frac{(\alpha \delta - \beta \gamma)(f_1g_2-f_2g_1)}{(\alpha \delta + \beta \gamma)(f_1g_2+f_2g_1)}\right|=\left|\frac{\alpha \delta - \beta \gamma}{\alpha \delta + \beta \gamma}\right|d(\bar{f},\bar{g}).\]

By Proposition \ref{prop:segments}, since the points $Ma$ and $Mb$, are coplanar with $\bar{f},\bar{g}$, the vector lines $(0,Ma)$ and $(0,Mb)$ intersect the segment $[a_1,b_1]$. Thus, we can apply Proposition \ref{prop:expression_distance} :
\[d(M\cdot \Pi(a), M\cdot \Pi(b))=\left|\frac{\alpha \delta - \beta \gamma}{\alpha \delta + \beta \gamma}\right|\]
Since $d(M\cdot \Pi(a), M\cdot \Pi(b))\leq c(M)$, we can finally conclude :
\[ d(M\cdot \bar{f},M \cdot \bar{g}) \leq c(M)d(\bar{f},\bar{g}).\]
\end{proof}
\begin{cor}
Consider a linear map $M:E\rightarrow E$ such that $M(C)\subset C$. Then the projective action of $M$ on $\Pi(C)$ is $c(M)$-Lipschitz with respect to the distance $d$. 
\end{cor}
We recall that a linear positive mapping $M:E\longrightarrow E$ is $A$-uniformly positive whenever there exists $h\in C$, such that for any $f\in C$, there exists $b(f)\in \RR_+$ such that :
\[ \frac{1}{A} b(f) h \leq Mf\leq A b(f) h.\]
Notice that the inequality remains valid replacing $A$ by $A'\geq A$. We note $A^*(M)=\inf\{A\geq 1 | $M$\text{ is } A-\text{uniformly positive}\}$.
Then the following characterisation of strictly contracting applications holds :
\begin{prop} \label{prop:strict contraction}
Consider a linear map $M:E\rightarrow E$, such that $M(C)\subset C$.
Then $c(M)<1$ if and only if $M$ is a uniformly positive linear mapping. Moreover when $c(M)<1$, $c(M)=\psi(A^*(M))$
where $\psi:s\in[1,\infty)\mapsto \frac{1-s^{-2}}{1+s^{-2}}=\phi(s^{-2})\in[0,1)$.
\end{prop}
\begin{proof}
First, assume that $M$ is $A$-uniformly positive for some $A$.
Let $f,g\in C$. Let $b(f), b(g)\in \RR$ so that
$$A^{-1}b(f)h \leq Mf \leq A b(f) h$$
$$A^{-1}b(g)h \leq Mg \leq A b(g) h.$$

Let $f'=Mf/ b(f)$ and $g'=Mg/b(g)$. Then \[ h\leq f',g'\leq A h,\]
hence
\[A^{-1}g'\leq h \leq Af' \text{ and } A^{-1}f'\leq h \leq Ag'.\]
Thus $\aleph(f',g'),\aleph(g',f')\geq A^{-1}$ and $m(f',g')\geq A^{-2}.$ Since $\phi:s\in [0,1]\mapsto\frac{1-s}{1+s}\in [0,1]$ is decreasing and one-to one, then $\psi$ is increasing and one-to-one. As a consequence,
\[ d(M\cdot \bar f, M\cdot \bar g)=\phi(m(f,g))=\phi(m(f',g'))\leq\phi(A^{-2})=\psi(A)<1.\]
Since $c(M)=\sup \left\{d(M\cdot \bar{f}, M\cdot \bar{g} ),\bar{f},\bar{g}\in \Pi(C)\right\},$ this yields $c(M)\leq \psi(A)<1.$ Since $\psi$ is continuous, this yields $c(M)\leq \psi(A^*(M))$.
\\ Conversely, suppose that $0\leq c(M)<1$. We want to show that $M$ is uniformly positive. Since $\psi$ is one-to-one from $[1,\infty)$ onto $[0,1)$, let us set $A=\psi^{-1}(c(M))$, thus $c(M)=\frac{1-A^{-2}}{1+A^{-2}}$.
Let $\bar{f},\bar{g}\in \Pi(C)$ and choose $f'\in M\cdot \bar{f}, g'\in M\cdot \bar{g}$, such that $(f',g')\cap C$ is a segment of line with endpoints $u,v$. We assume that $f'=s u +(1-s)v$ and $g'=t u +(1-t) v$ with $0\leq s \leq t \leq 1$. Since we have assumed that $c(M)<1$, then $d(M\cdot \bar{f},M\cdot \bar{g})\leq c(M)< 1$. Therefore, $f'$ and $g'$ are distinct from $u,v$ (ie $0< s,t< 1$) and
$$d(M\cdot \bar{f},M\cdot \bar{g})=\frac{1-m(f',g')}{1+m(f',g')}\leq c(M)=\frac{1-A^{-2}}{1+A^{-2}}.$$
Since $s\mapsto \frac{1-s}{1+s}$ is strictly decreasing on $(0,1)$, this yields
\[m(f',g')\geq \frac{1}{A^2}. \]
Let us choose $r>0$ such that $\aleph(r f',g')=r \aleph(f',g')= A^{-1}$, and set $f''=r f', g''=g'$. Then it holds,
\[ A^{-2}\leq m(f',g')=m(f'',g'')=\aleph(f'',g'')\aleph(g'',f'')=A^{-1}\aleph(g'',f'').\]
Thus $\aleph(f'',g''),\aleph(g'',f'') \geq A^{-1}$, which yields $f''\leq Ag''$, and $g''\leq Af''$.
In other words : if $c(M)<1$, then there exists $A>0$ such that, for any $\bar{f},\bar{g}\in \Pi(C)$, there exists $f''\in M\cdot \bar{f}, g''\in M\cdot \bar{g}$, satisfying $A^{-1}g''\leq f'' \leq A g''$.
Let us now consider an arbitrary function $g\in C$, and fix $h=Mg$.
For each $f\in C$, there exists two positive reals $\alpha(f),\beta(f)>0$ such that $f''=\alpha(f) Mf $ and $g''=\beta(f) Mg=\beta(f)h$ satisfy \[ A^{-1} g'' \leq f'' \leq A g''.\]
This yields :
\[ A^{-1} \frac{\beta(f)}{\alpha(f)} f \leq Mf \leq A \frac{\beta(f)}{\alpha(f)} h.\]
In other words, setting $b(f)=\frac{\beta(f)}{\alpha(f)}$,
\[ {A}^{-1} b(f) h \leq Mf \leq A b(f) h.\]
Thus, $M$ is $A$-uniformly positive, and $A^*(M)\leq A=\psi^{-1}(c(M))$, thus $\psi(A^*(M))\leq c(M)$.

\end{proof}
\section{Examples}
\subsection{The finite dimensional case}
Let us focus in this subsection on the case of a finite set $\XX=\{1,\dots, d\}$. We set $E=\RR^d$, $C=(\RR_+)^d-\{0\}$, and consider $d\times d$ matrices with non-negative coefficients, acting linearly on row vectors. In this case the definition of uniform positivity is equivalent to considerations on the zero-coefficients of $M$ :
\begin{prop}
Let $M=(M_{ij})_{1\leq i,j\leq d}$ be a $d\times d$ matrix with non-negative coefficients. Then $M$ is uniformly positive if and only if for each $i,j\in \{ 1,\dots d\}$ such that $M_{ij}=0$, either the whole $i-$th line or the whole $j-$th row is zero.
\end{prop}
\begin{proof}
Let us assume that $M$ is uniformly positive, then there exists $A>0$, and two functions $x\mapsto h(x), y\mapsto b(y)$ satisfying for every row vector $y$ and line vector $x$
\begin{equation} \label{eq:unifpositivity}A^{-1}h(x) b(y)\leq xMy \leq Ah(x) b(y).\end{equation}
Thus if $M_{i,j}=e_iMe_j=0$, either $h(e_i)=0$ or $ b(e_j)=0$. This implies respectively that the $j$-th row or the $i$-th column of $M$ is $0$.
\\ Conversely, if we know that $M$ has no zero coefficients except in whole zero line or rows, then letting $i_0,j_0$ be the respective indexes of a non-zero line and a non-zero row, then setting $h(e_i)=M_{ij_0}$ and $ b(e_j)=M_{i_0j}$,and
\[ 0<A=\max\left( \max_{i,j} \frac{M_{ij}}{M_{ij_0}M_{i_0j}}, \max_{i,j}\frac{M_{ij_0}M_{i_0j}}{M_{ij}}\right)<\infty,\]
we prove that \eqref{eq:unifpositivity} is satisfied for any vectors $x,y$ of the canonical basis. It can be extended to the whole positive cone by linearity.
\end{proof}
This provides a natural necessary and sufficient condition for a matrix to be a strictly contracting projective operator on the projective cone $\Pi(C)$.
\begin{cor}\label{cor:cns_contraction_matrix}
Let $M$ be a $d\times d$ matrix with non negative coefficients. Then
\begin{itemize}
    \item[i)] $M(C)\subset C$ if and only if there is a non zero coefficient on each row of $M$
    \item[ii)] Suppose that $M(C)\subset C$. Then $c(M)<1$ if and only if for each $i,j$ such that $M_{ij}=0$, the whole $j-$th column of $M$ is $0$.
\end{itemize}
\end{cor}

Additionally, the following Proposition, stated and proved in \cite[Lemma 10.7]{hennion_limit_1997} claims that the constant $c(M)$ can be explicitly computed in terms of the coefficients of $M$.
\begin{prop}
    Let $M$ is be a $d\times d$ matrix, with non positive coefficients, and no zero-row. Then 
\[ c(M)=\sup\left\{d(M\cdot x,M\cdot y), x,y\in \Pi(C)\right\}=\max \left\{d(M\cdot e_i,M\cdot e_j),i,j\in \left\{1,\dots n\right\}\right\},\]
where $(e_i)_{1\leq d}$ are the unit vectors of the canonical basis of $\RR^d.$ 
As a consequence, the constant $c(M)$ can be explicitly expressed in terms of the coefficients of $M=\left(M_{i,j}\right)_{1\leq i,j \leq d}$ as follows :
\begin{equation}\label{eq:formule coeff} c(M)=\underset{1\leq i,j,k,l\leq d}{\max} \frac{\left \vert M_{ki}M_{lj}-M_{kj}M_{li}\right\vert}{M_{ki}M_{lj}+M_{kj}M_{li}}.\end{equation}
with the convention $0/0=1.$
\end{prop}
One can easily check that the necessary and sufficient condition for $c(M)<1$ provided in Corollary \ref{cor:cns_contraction_matrix} is consistent with the expression of $c(M)$ in \eqref{eq:formule coeff}.

\subsection{The density case}
Let $\XX=[0,1]$ and let $E$ be the set of continuous real functions on $\XX$. Let us consider linear operators of the form \[M_K f(x)=\int K(x,y)f(y)dy,\]
where $(x,y)\mapsto K(x,y)$ is a continuous map from $[0,1]$ to $\RR_+$. Under these assumptions $M_K$ is a bounded linear operator on $E$.
Assume additionally that 
\begin{itemize}
    \item For each $x,y\in\RR$, $K(x,y)\geq 0$
    \item For each $y\in\RR$, there exists $x\in\RR$ such that $K(x,y)>0$.
\end{itemize}
These conditions ensure that $M_K$ is a positive operator on $E$, satisfying $M_K(C)\subset C$.
\begin{prop}
Under these assumptions $c(M_K)<1$ if and only if there exists $g_1,g_2\in \RR_+^\XX, A\in \RR_+^*$ such that
\begin{equation}\label{eq:unifpos noyau} A^{-1} g_1(x)g_2(y)\leq K(x,y)\leq A g_1(x)g_2(y).\end{equation}
When \eqref{eq:unifpos noyau} holds for some $A$, noting $A^*$ the infimum of the values of $A$ satisfying \eqref{eq:unifpos noyau}, then it holds 
\[ c(M_K)=\psi^{-1}(A^*).\]
\end{prop}
By continuity and compacity, it is sufficient to have $K(x,y)>0$ for all $x,y$ in order for equation \eqref{eq:unifpos noyau} to hold. However, it is clearly not necessary, since any product of the form $K(x,y)=g_1(x)g_2(y)g_3(x,y)$ with $g_1,g_2,g_3$ continuous and non negative, and $g_3$ positive on $[0,1]^2$ satisfies \eqref{eq:unifpos noyau} and is therefore uniformly positive. 

\begin{proof}
If there exists $A,g_1,g_2$ satisfying equation \eqref{eq:unifpos noyau}, then for any $x\in \XX$, any $f\in E$
\[A^{-1} g_1(x) \int g_2(y) f(y)dy \leq M_K f (x) \leq A g_1(x) \int g_2(y) f(y)dy, \]
thus $M_K$ is $A$-uniformly positive with $h=g_1$ and $ b(f)=\int g_2(y) f(y)dy$.
\\ Conversely, if $M_K$ is $A$-uniformly positive, let $h\in \RR^\XX$ and for each $f\in E$, $ b(f)>0$ such that
 \[ A^{-1} b(f) h \leq Mf \leq  
 K b(f) h.\]
 Let us note, for $\varepsilon>0$ and $y\in \XX$,
 \[g_{\varepsilon,y}=\frac{\mathds{1}_{[y-\varepsilon,y+\varepsilon]\cap \XX}}{\int \mathds{1}_{[y-\varepsilon,y+\varepsilon]\cap \XX}(z)dz}.\]
 Then it holds, for any $x,y\in\XX$, by continuity of $y\mapsto K(x,y)$
 \[ K(x,y)=\lim_{\varepsilon\rightarrow 0} M_K(g_{\varepsilon,y})(x).\]
 The uniform positivity assumption yields however that 
 \[ A^{-1} M_K(g_{\varepsilon,y})(x) \leq h(x) b(g_{\varepsilon,y})\leq A M_K(g_{\varepsilon,y})(x).\]
 Therefore 
 \[ A^{-1} K(x,y) \leq h(x) \underset{\varepsilon \rightarrow 0}\limsup \, b(g_{\varepsilon,y}) \leq A K(x,y) .\]
 Setting $g_1=h$ and $g_2(y)=\underset{\varepsilon \rightarrow 0}\limsup \,b(g_{\varepsilon,y})$, we get 
 \[A^{-1}g_1(x)g_2(y)\leq K(x,y) \leq A g_1(x)g_2(y).\]
This proves that $M_K$ is $A$ uniformly positive if and only if \eqref{eq:unifpos noyau} holds for some functions $g_1,g_2$. Proposition \ref{prop:strict contraction} allows to conclude the proof.
 \end{proof}

\section{Acknowledgements}
I have received support from the Chair "Modélisation Mathématique et Biodiversité" of VEOLIA, Ecole Polytechnique, MnHn, FX,  as well as from the ANR project NOLO (ANR 20-CE40-0015), funded by the French ministry of research. I would like to warmly thank my PhD supervisors Vincent Bansaye and Marc Peigné for their support and feedbacks on this note.

\printbibliography

\end{document}